\documentclass[12pt]{amsart}
\usepackage{amsthm, amssymb}
\usepackage{amsfonts, amscd}
\usepackage{epsfig,multicol}
\theoremstyle{plain} \numberwithin{equation}{section}

\theoremstyle{definition}

\begin{document}
\title[Corrections to: Involutions fixing $\Bbb{RP}^{\text{odd}}\sqcup P(h,i)$, II]
{\large \bf Corrections to: Involutions fixing
$\Bbb{RP}^{\text{odd}}\sqcup P(h,i)$, II}
\author[Bo Chen and Zhi L\"u]{Bo Chen and Zhi L\"u}
\footnote[0]{ 
\endgraf
\endgraf
  Supported by  grants from NSFC (No. 10371020). }
\address{Institute of Mathematics, School of Mathematical Sciences, Fudan University, Shanghai,
200433, P.R. China.} \email{042018018@fudan.edu.cn}

\address{Institute of Mathematics, School of Mathematical Sciences, Fudan University, Shanghai,
200433, P.R. China.} \email{zlu@fudan.edu.cn}
\date{}
\maketitle

The purpose of this note is to correct statements of some
assertions in \cite{l}. The mistake occurs in the argument of the
case in which the normal bundle $\nu^k$ over $P(h,i)$ is
nonstandard. Specifically, some incorrect calculations first
happen in the arguments of the cases $u=0$ and $u>1$ of page 1309
(in the proof of Lemma 3.4 of \cite{l}). This leads to the loss of
the existence of some involutions with nonstandard normal bundle
$\nu^k$ in those two cases, so that the statements of Lemma 3.4
and Proposition 3.4 are incorrect, and so is  part of the
statement of Theorem 2.3 in \cite{l}.

    \vskip .2cm

Following the notations of \cite{l},    Lemma 3.4 in \cite{l}
should be corrected as follows.

    \vskip .2cm
\noindent {\bf Lemma.} {\em If $\nu^k$ is nonstandard, then $h=2$
with the following possible cases:

$(A)$ for $u=0$, one has that $(k,a)=(2,2)$,
$\nu^2=\tau\otimes\eta$ and $i\equiv 3\mod 4$;

$(B)$ for $u=1$, one has that $(k,a)=(4,1)$ and
$\nu^4=\tau\oplus(\tau\otimes\eta)$;

$(C)$ for $u>1$, one has that either $a=1$ and $6\leq k\ (\equiv
2\mod 4)$ with $\nu^k$ stably cobordant to
$3\xi\oplus(\tau\otimes\eta)$ or $a=3$ and $4\leq k\ (\equiv 0\mod
4)$ with $\nu^k$ stably cobordant to $\xi\oplus(\tau\otimes\eta)$,
where $\xi$ is a 1-plane bundle  over $P(2,2^u(2v+1))$, $\eta$ is
a 2-plane bundle  over $P(2,2^u(2v+1))$, and $\tau$ is the 2-plane
bundle (the tangent bundle of $\Bbb{RP}^2$ pulled back to
$P(2,2^u(2v+1))$).}

\vskip .2cm

{\em Note.} Stong in \cite{s} found the strange tensor product
$\tau\otimes\eta$ over the Dold manifold with the total class
$1+c+c^2+d$.

\begin{proof}
Since the mistake in the proof of Lemma 3.4 of \cite{l} only
occurs in the cases $u=0$ and $u>1$ of page 1309 but other
arguments are true, one needs to merely show (A) and (C).

\vskip .2cm If $u=0$, then $a$ is even and $k=2$ by Lemma 3.1 in
\cite{l}. To ensure $k=2$, from (3.2) in \cite{l} one must have
$a=2$, so the total class $w(\nu^2)=1+c+c^2+d$. Thus,
$\nu^2=\tau\otimes\eta$. By direct computation, one has that
$$
w[0]_1=\begin{cases}
c & \text{ on } P(2,i)\\
\alpha  & \text{ on } \Bbb{RP}^2
\end{cases}
\text{ and }  w[0]_2=\begin{cases}
ce+c^2+d+{{i+3}\choose 2}c^2 & \text{ on } P(2,i)\\
\alpha^2+e^2  & \text{ on } \Bbb{RP}^2.
\end{cases}
$$
Form the class $$\hat{w}_2=w[0]_2+e^2+w[0]_1^2=\begin{cases}
ce+e^2+d+{{i+3}\choose 2}c^2 & \text{ on } P(2,i)\\
0 & \text{ on } \Bbb{RP}^2
\end{cases}$$
one has that the value of $\hat{w}_2e^{2i+1}$ on $\Bbb{RP}^2$ is
zero, so the value of this on $P(2,i)$ is zero, too. Thus
$$0=\hat{w}_2e^{2i+1}[\Bbb{RP}(\nu^2)]={{1+c+d+{{i+3}\choose 2}c^2}\over{1+c+c^2+d}}[P(2,i)]=
1+{{i+3}\choose 2}$$ and so $i\equiv 3\mod 4$.

\vskip .2cm

If $u>1$,  then, by Lemma 3.1 in \cite{l}, $a$ is odd  and $k$ is
even. Further, $a=1$ or 3 since $h=2$. Now by direct calculations,
\begin{eqnarray*}
{1\over {w(\nu^k)}}[P(2,
2^u(2v+1))]&=&{1\over{(1+c)^a(1+c+d)(1+{{c^2d}\over{1+d}})}}[P(2,
2^u(2v+1))]\\
&=& {{d^{2^u(2v+1)}}\over{(1+c)^{a+1}}}[P(2,
2^u(2v+1))]\\
&=&\begin{cases} 1 & \text{if $a=1$}\\
0 & \text{if $a=3$}
\end{cases}
\end{eqnarray*}
and
$${1\over{w(\nu^{j+1})}}[\Bbb{RP}^2]={1\over{(1+\alpha)^k}}[\Bbb{RP}^2]=\begin{cases} 1 & \text{if $k\equiv 2\mod 4$}\\
0 & \text{if $k\equiv 0\mod 4$}.
\end{cases}$$
Thus, $a=1$ if and only if $k\equiv 2\mod 4$, and $a=3$ if and
only if $k\equiv 0\mod 4$. When $a=1$, one has
$w(\nu^k)=(1+c)^3(1+c+c^2+d)$ so $k>4$ and $\nu^k$ is stably
cobordant to $3\xi\oplus(\tau\otimes\eta)$; when $a=3$, one has
$w(\nu^k)=(1+c)(1+c+c^2+d)$ so $k>2$ and $\nu^k$ is stably
cobordant to $\xi\oplus(\tau\otimes\eta)$.
\end{proof}

Next, Proposition 3.4 in \cite{l} should be corrected as follows.

\vskip .2cm

\noindent {\bf Proposition.} {\em The involution
$(M^{2^{u+1}(2v+1)+k+h},T)$ fixing
$\Bbb{RP}^{2^{u+1}(2v+1)+k-1}\sqcup P(h,2^u(2v+1))$ with $\nu^k$
nonstandard  exists only for the following four cases:

$(i)$ $(h, u,k,a)=(2,0,2,2)$, $\nu^2=\tau\otimes\eta$ and $v$ is
odd;

$(ii)$ $(h, u,k,a)=(2,1,4,1)$ and
$\nu^2=\tau\oplus(\tau\otimes\eta)$;

$(iii)$ $(h,a)=(2, 1)$ with $u>1$,  $k\equiv 2\mod 4$ is in the
range $6\leq k\leq Y_1$ and $\nu^k$ is stably cobordant to
$3\xi\oplus(\tau\otimes\eta)$, where $Y_1\leq 2^{u+1}-2$;

$(iv)$ $(h,a)=(2, 3)$ with $u>1$,  $k\equiv 0\mod 4$ is in the
range $4\leq k\leq Y_2$ and $\nu^k$ is stably cobordant to
$\xi\oplus(\tau\otimes\eta)$, where $Y_2\leq 2^{u+1}$.}

\vskip .2cm

 {\em Note.} Proposition 3.4 in \cite{l} only indicates the existence of the
involution of  case (ii) in the above proposition, and its proof
is true. However, as stated in the above proposition, actually
there are also other cases in which the involutions with $\nu^k$
nonstandard  exist.

\begin{proof} First, by the above lemma, one has $h=2$.
As stated in the introduction of \cite{l}, it suffices to discuss
the existence of involutions
$(\bar{M}^{2^{u+1}(2v+1)+k+2},\bar{T})$ fixing $\Bbb{RP}^2$ with
normal bundle $\nu^{2^{u+1}(2v+1)+k}$ having
$w(\nu^{2^{u+1}(2v+1)+k})=(1+\alpha)^{2^{u+1}+k}$ and
$P(h,2^u(2v+1))$ with $\nu^k$ nonstandard. In a similar way to the
argument of  case (ii) as shown in the proof of Proposition 3.4 of
\cite{l}, one can easily prove that the  involution with $\nu^k$
nonstandard exists for the following cases:

\vskip .2cm

(a) $(h, u,k,a)=(2,0,2,2)$, $\nu^2=\tau\otimes\eta$ and $v$ is
odd;

(b) $(h, k, a)=(2,6,1)$ with $u>1$ and
$w(\nu^6)=(1+c)^3(1+c+c^2+d)$;

(c) $(h, k, a)=(2,4,3)$ with $u>1$ and
$w(\nu^4)=(1+c)(1+c+c^2+d)$,

\vskip .2cm

\noindent which means that the above proposition holds for  case
(i),  case (iii) with $k=6$, and  case (iv) with $k=4$. In
particular,  the same argument as above can also show that
$\bar{M}^{2^{u+1}(2v+1)+k+2}$ is cobordant to zero in  case (b)
with $u>2$ and  case (c). Furthermore, one can apply the
$\Gamma$-operation to $(\bar{M}^{2^{u+1}(2v+1)+k+2},\bar{T})$ to
obtain more involutions with nonstandard $\nu^k$. Thus, it remains
to estimate the upper bound of $k$ in  cases (iii) and (iv). If
$u>1$, by direct computations, one has that
$$w[0]_4=\begin{cases}
c^2d+cde+de^2+d^2+{{a+1}\choose 2}c^2e^2 & \text{ on } P(2,
2^u(2v+1))\\
{{a+1}\choose 2}\alpha^2e^2 & \text{ on } \Bbb{RP}^2.
\end{cases}$$
Form the class
$$\hat{w}_4=w[0]_4+{{a+1}\choose 2}w[0]_1^2e^2=\begin{cases}
c^2d+cde+de^2+d^2 & \text{ on } P(2,2^u(2v+1))\\
0 & \text{ on } \Bbb{RP}^2.
\end{cases} $$
For  case (iii), if $k>2^{u+1}-2$, one has that the value of
$$\hat{w}_4^{2^u(v+1)}e^{1+k-2^{u+1}}$$
on $\Bbb{RP}^2$ is zero, but the value of this on $P(2,
2^u(2v+1))$ is
\begin{eqnarray*}
\hat{w}_4^{2^u(v+1)}e^{1+k-2^{u+1}}[\Bbb{RP}(\nu^k)]&=&
{{d^{2^u(v+1)}(1+c+c^2+d)^{2^u(v+1)}}\over {w(\nu^k)}}[P(2,
2^u(2v+1))]\\
&=& {{d^{2^u(v+1)}(1+c+c^2+d)^{2^u(v+1)}}\over
{(1+c)^3(1+c+c^2+d)}}[P(2,2^u(2v+1))]\\
&=& (1+c)d^{2^u(v+1)}(1+c+c^2+d)^{2^u(v+1)-1}[P(2,2^u(2v+1))]\\
&=& d^{2^u(2v+1)}{{2^u(v+1)-1}\choose
{2^uv}}(1+c)^{3(2^u-1)+1}[P(2,2^u(2v+1))]\\
&=& c^2d^{2^u(2v+1)}[P(2,2^u(2v+1))]\\
&=& 1
\end{eqnarray*}
which leads to a contradiction. Thus, one has that $k\leq
2^{u+1}-2$ so $Y_1\leq 2^{u+1}-2$.

\vskip .2cm

For  case (iv), if $k>2^{u+1}$, one has that the value of
$$\hat{w}_4^{2^u(v+1)}(1+w[0]_1)^2e^{k-2^{u+1}-1}$$
on $\Bbb{RP}^2$ is zero, but the value of this on $P(2,
2^u(2v+1))$ is
\begin{eqnarray*}
\hat{w}_4^{2^u(v+1)}e^{1+k-2^{u+1}}[\Bbb{RP}(\nu^k)]&=&
{{d^{2^u(v+1)}(1+c+c^2+d)^{2^u(v+1)}(1+c)^2}\over {w(\nu^k)}}[P(2,
2^u(2v+1))]\\
&=& {{d^{2^u(v+1)}(1+c+c^2+d)^{2^u(v+1)}(1+c)^2}\over
{(1+c)(1+c+c^2+d)}}[P(2,2^u(2v+1))]\\
&=& (1+c)d^{2^u(v+1)}(1+c+c^2+d)^{2^u(v+1)-1}[P(2,2^u(2v+1))]\\
&=& 1.
\end{eqnarray*}
This is impossible. Thus, one has that $k\leq 2^{u+1}$ so $Y_2\leq
2^{u+1}$.
\end{proof}

Finally, combining the above lemma and proposition, the correct
statement of Theorem 2.3 in \cite{l} should be the following.

\vskip .2cm

\noindent {\bf Theorem.} {\em Suppose that $(M^{j+q},T)$ fixes
$\Bbb{RP}^j$ with normal bundle $\nu^q$ having
$w(\nu^q)=(1+\alpha)^q$ with odd $q>1$, and $P(h,i)$ with normal
bundle $\nu^k$ having
$w(\nu^k)=(1+c)^a(1+c+d)^bw(\rho)^\varepsilon$. Let $2^A\leq h\leq
2^{A+1}$ and write $i=2^u(2v+1)$. Then $(b, q, j)=(1, h+1,
2i+k-1)$,  $k$ is even with $2\leq k\leq \begin{cases} 2^{u+1}+2 & \text{ if } u=1\\
2^{u+1} & \text{ if } u\not=1
\end{cases}$ and  $i+a$ is odd.

\vskip .2cm

{\em (I)} When $\varepsilon=0$ (i.e., $\nu^k$ is standard), one
has that

$(1)$  $a<2^u$.

$(2)$ $j+1\equiv i+a+1\mod 2^{A+1}$ and $i+k\equiv a+1\mod
2^{A+1}$. In particular,

\ \ $(a)$ for $u\leq A$, $k=2^u+a+1$ and $2^{u+1}(v+1)\equiv 0\mod
2^{A+1}$;

\ \ $(b)$ for $u>A$, $k\equiv a+1\mod 2^{A+1}$.

\noindent Further, $(M^{j+q},T)$ with standard $\nu^k$ exists for
$k$ in a range $X_1\leq k\leq X_2$, and is cobordant to
$$\Gamma^{k-2a-2}(P(h,N^{i+a+1}), T_{N^{i+a+1}})\sqcup
(\Bbb{RP}^{j+h+1}, T_{h+1})$$ where $2\leq X_1,X_2\leq
2k_0=\begin{cases} 2^{u+1}+2 & \text{ if } u=1\\
2^{u+1} & \text{ if } u\not=1
\end{cases}$ and more precisely

\ \ $(c)$ for $u\leq A$, $X_1=a+2$ and $X_2\leq 2^u+a+1$;

\ \ $(d)$ for $u>A$, $2\leq X_1\leq h+2$ and $X_2\leq
2^{u+1}-(h-\text{common}(h,a))$ where common$(h,a)$ is the common
part of the 2-adic expansions of $h$ and $a$.

\vskip .2cm

{\em (II)} When $\varepsilon\not=0$ (i.e., $\nu^k$ is
nonstandard), one has $h=2$. Further, $(M^{j+q},T)$ with
nonstandard $\nu^k$ exists only for the following cases:

$(1)$ $(u,k,a)=(0,2,2)$, $\nu^2=\tau\otimes\eta$ and $v$ is odd;

$(2)$ $(u,k,a)=(1,4,1)$ and $\nu^2=\tau\oplus(\tau\otimes\eta)$;

$(3)$ $a=1$ with $u>1$,  $k\equiv 2\mod 4$ is in the range $6\leq
k\leq Y_1$ and $\nu^k$ is stably cobordant to
$3\xi\oplus(\tau\otimes\eta)$, where $Y_1\leq 2^{u+1}-2$;

$(4)$ $a=3$ with $u>1$,  $k\equiv 0\mod 4$ is in the range $4\leq
k\leq Y_2$ and $\nu^k$ is stably cobordant to
$\xi\oplus(\tau\otimes\eta)$, where $Y_2\leq 2^{u+1}$. }

\vskip .2cm

In concluding this note, it should be pointed out that there is an
additional number 384 in line 18 of  page 4555 in \cite{l*}, which
should be omitted.

\end{document}